\documentclass[11pt,a4paper]  
 {article}
\usepackage{amssymb,latexsym, amsmath, bm, amsfonts, tikz,fancyhdr,sectsty}
\usepackage{pgfplots}
\usetikzlibrary{patterns}

\usepackage{lettrine}
\newtheorem{theorem}{Theorem}[section]
\newtheorem{lemma}[theorem]{Lemma}

\newtheorem{corollary}[theorem]{Corollary}

 \title{\bf Generalised differences and  multiplier 
   operators  in  $\boldsymbol{L^2({\mathbb R})}$}
\author{Rodney Nillsen \\
\small{School of Mathematics and Applied Statistics}\\
\small{University of Wollongong,
  Northfields Avenue, Wollongong}\\
\small{New South Wales, 2522, AUSTRALIA}\\
\small{email address: nillsen@uow.edu.au}\\}
 \date{\ }
  
\begin{document}
\maketitle
 \begin{abstract}
\noindent Let $\alpha, \beta\in {\mathbb R}$ be given  and let $s\in {\mathbb N}$ also be given.  Let $\delta_x$ denote the Dirac measure at $x\in {\mathbb R}$, and  let $\ast$ denote convolution. If $f\in L^2({\mathbb R})$, and if  there are  $u\in {\mathbb R}$  and   $g\in L^2({\mathbb R})$ such that
 \[ f=\left[\left([e^{iu\left(\frac{\alpha-\beta}{2}\right)}+e^{-iu\left(\frac{\alpha-\beta}{2}\right) } \right)\delta_0-  \left(\  e^{iu\left(\frac{\alpha+\beta}{2}\right)}\,\delta_{ u}+e^{-iu\left(\frac{\alpha+\beta}{2}\right)}\,\delta_{-u}\right)\right]^s\ast g,\]
 then $f$ is called a \emph{generalised $(\alpha,\beta)$-difference of order $2s$}, or simply a \emph{generalised  difference}.
     We denote by ${\mathcal D}_{\alpha,\beta,s}({\mathbb R})$ the vector space of all functions $f$ in  $L^2({\mathbb R})$ such that   $f$ is a finite sum of  generalised  $(\alpha,\beta)$-differences  of order $2s$.  
 It is shown that every function in  ${\mathcal D}_{\alpha,\beta,s}({\mathbb R})$ is a sum of $4s+1$  generalised $(\alpha,\beta)$-differences of order $2s$.  If we let ${\widehat f}$ denote the Fourier transform of a function $f\in L^2({\mathbb R})$,
 then  ${\mathcal D}_{\alpha,\beta,s}({\mathbb R})$ is a weighted $L^2$-space under the Fourier transform, and its  inner product $\langle\, \,,\,\rangle_{\alpha,\beta,s}$ is given by 
 \[{\langle f,g\rangle}_{\alpha,\beta,s} =\int_{-\infty}^{\infty}\left(1+\frac{1}{(x-\alpha)^{2s}(x-\beta)^{2s}}\right)  {\widehat f}(x){\overline{{\widehat g}(x)}}dx.\]
 Letting $D$ denote  differentiation and letting $I$ denote the identity operator,  the operator $(D^2-i(\alpha+\beta)D-\alpha\beta I)^s$ is bounded and invertible, mapping the Sobolev subspace of order $2s$ of $L^2({\mathbb R})$ onto  the Hilbert space ${\mathcal D}_{\alpha,\beta,s}({\mathbb R})$.
  \end{abstract}
\let\thefootnote\relax
 
  \section{Introduction}
\setcounter{equation}{0}
Let ${\mathbb R}$ denote the set of real numbers, let ${\mathbb T}$ denote the set of complex numbers of modulus $1$, and let $G$ denote either ${\mathbb R}$ or ${\mathbb T}$. Note that in some contexts ${\mathbb T}$ may be identified with the interval $[0,2\pi)$ under the mapping $t\longmapsto e^{it}$ (some comments on this are in \cite[page 1034]{ross1}). Then $G$ is a group and its identity element we denote by $e$, so that $e=0$ when $G={\mathbb R}$ and $e=1$ when $G={\mathbb T}$. Let  ${\mathbb N}$ denote the set of natural numbers,  ${\mathbb Z}$  the set of integers, and let $s\in {\mathbb N}$.  The Fourier transform of $f\in L^2(G)$  is denoted by ${\widehat f}$, and is given by ${\widehat f}(n)=(2\pi)^{-1}\int_0^{2\pi}f(e^{it})e^{-int}$ for $n\in {\mathbb Z}$ (in the case of $\mathbb T$), and by  the extension to all of $L^2(\hbox{$\mathbb R$})$ of the transform given by ${\widehat f}(x)=\int_{-\infty}^{\infty}e^{-ixu}f(u)\,du$ for $x\in {\mathbb R}$ (in the case of ${\mathbb R}$). Let $M(G)$ denote the family of bounded Borel measures  on $G$. If $x\in G$ let $\delta_x$ denote the Dirac measure at $x$, and let $\ast$ denote convolution in $M(G)$.

We call a function $f\in L^2(G)$ a \emph{difference of order} $s$ if there is a function $g\in L^2(G)$ and $u\in G$ such that $f=(\delta_e-\delta_u)^s\ast g$. The functions in $L^2(G)$ that are a sum of a finite number of differences of order $s$ we denote  by ${\mathcal D}_s(G)$. Note that ${\mathcal D}_s(G)$ is a vector subspace of $L^2(G)$.  Now in the case of ${\mathbb T}$ it was shown by Meisters  and Schmidt \cite{meisters1} that
\begin{equation}{\mathcal D}_1({\mathbb T})=\Bigl\{f:f\in L^2({\mathbb T})\ {\rm and}\ {\widehat f}(0)=0\Bigr\}, \nonumber
\end{equation}
and that every function in ${\mathcal D}_1({\mathbb T})$ is a sum of $3$ differences of order $s$.  It was shown in \cite{nillsen1} that, for all $s\in {\mathbb N}$,
\begin{equation}{\mathcal D}_s({\mathbb T})={\mathcal D}_1({\mathbb T})=\Bigl\{f:f\in L^2({\mathbb T})\ {\rm and}\ {\widehat f}(0)=0\Bigr\},\label{eq:characterisation1}\end{equation}
and that every function in ${\mathcal D}_s({\mathbb T})$ is a sum of $2s+1$ differences of order $s$. It was also shown in \cite{nillsen1} that
\begin{equation}{\mathcal D}_s({\mathbb R})=\left\{f:f\in L^2({\mathbb R})\ {\rm and}\  \int_{-\infty}^{\infty}\frac{|{\widehat f}(x)|^2}{|x|^{2s}}\,dx<\infty\right\},
\label{eq:characterisation2}
\end{equation}
and  again, that every  function in ${\mathcal D}_s({\mathbb R})$ is a sum of $2s+1$ differences of order $s$. Further   results  related to the work of Meisters and Schmidt in \cite{meisters1} may be found in \cite{bourgain1, johnson1, meisters2, nillsen2}.

The Sobolev space of order $s$  in $L^2(G)$ is the space of all functions $f\in L^2(G)$ such that $D^s(f)\in L^2(G)$, where $D$ denotes differentiation in the sense of Schwartz distributions. Then, $D^s$ is a \emph{multiplier operator} on $W^s({\mathbb T})$ with \emph{multiplier} $(in)^s$, in the sense that $D^s(f)^ {\widehat {\  }}(n)=(in)^s{\widehat f}(n),$ for all $f\in W^s({\mathbb T})$  and $n\in {\mathbb Z}$. Also, $D^s$ is a multiplier operator on $W^s({\mathbb R})$ with multiplier $(ix)^s$, the sense that $D^s(f)^{\widehat {\ }}(x)=(ix)^s{\widehat f}(x),$ for all $f\in W^s({\mathbb R})$  and $x\in {\mathbb R}$. Note   that $W^s({\mathbb T})$ is a Hilbert space whose norm $||\cdot ||_{{\mathbb T},s}$ derives from the inner product $\langle\,,\,\rangle_{{\mathbb T},s}$  where 
\[\langle f,g\rangle_{{\mathbb T},s}=\sum_{n=-\infty}^{\infty}\left(1+|n|^{2s}\right){\widehat f}(n){\overline{{\widehat g}(n)}}\, dx.\]Note also that $W^s({\mathbb R})$ is a Hilbert space whose norm $||\cdot ||_{{\mathbb R},s}$ derives from the inner product $\langle\,,\, \rangle_{{\mathbb R},s}$  where 
\[\langle f,g\rangle_{{\mathbb R},s}=\int_{-\infty}^{\infty}\left(1+|x|^{2s}\right){\widehat f}(x){\overline{{\widehat g}(x)}}\, dx.\]
Using these observations, together with Plancherel's Theorem, it is  easy to  verify  that
\begin{align}D^s(W^s({\mathbb T}))&=\Bigl\{f:f\in L^2({\mathbb T})\ {\rm and}\  {\widehat f}(0)=0\Bigr\},\ {\rm and\ that}\label{eq:ranges1}\\
D^s(W^s({\mathbb R}))
&= \left\{f:f\in L^2({\mathbb R})\ {\rm and}\  \int_{-\infty}^{\infty}\frac{|{\widehat f}(x)|^2}{|x|^{2s}}\,dx<\infty\right\}.\label{eq:ranges2}
\end{align}

In view of (\ref{eq:ranges1}) and (\ref{eq:ranges2}),  (\ref{eq:characterisation1}) together with   (\ref{eq:characterisation2})   can be regarded as describing the ranges of $D^s$ upon $W^s({\mathbb T})$ and $W^s({\mathbb R})$  as spaces  consisting of finite sums of differences of order $s$.  Corresponding   results have been obtained in \cite{nillsen3} for  operators  $(D^2-i(\alpha+\beta)D-\alpha\beta I)^s$ acting on $W^{2s}({\mathbb T})$, where $\alpha,\beta\in {\mathbb Z}$ and  $I$ denotes the identity operator. In this paper,  the  main aim  is to derive corresponding  results for the operator $(D^2-i(\alpha+\beta)-\alpha\beta I)^s$, where $\alpha,\beta\in {\mathbb R}$,   for the non-compact case of ${\mathbb R}$ in place of the compact group ${\mathbb T}$. 
Note that, in general, the range of a multiplier operator depends upon the behaviour of Fourier transforms at or around the zeros of the multiplier of the operator, as in (\ref{eq:ranges1}) and (\ref{eq:ranges2}). Note also that  $(D^2-i(\alpha+\beta)D-\alpha\beta I)^s$  is a multiplier operator whose multiplier is, in the case of ${\mathbb R}$,  $ -(x-\alpha)(x-\beta)$ with zeros at $\alpha$ and $\beta$.

Given $\alpha,\beta\in {\mathbb R}$ and $s\in {\mathbb N}$,  a \emph{generalised $(\alpha,\beta)$-difference of order $2s$}  is a function $f\in L^2(G)$  such that for some $g\in L^2(G)$ and $u\in G$  we have
 \begin{equation}
 f=\left[\left([e^{iu\left(\frac{\alpha-\beta}{2}\right)}+e^{-iu\left(\frac{\alpha-\beta}{2}\right) } \right)\delta_0-  \left(\  e^{iu\left(\frac{\alpha+\beta}{2}\right)}\,\delta_{u}+e^{-iu\left(\frac{\alpha+\beta}{2}\right)}\,\delta_{-u}\right)\right]^s\ast g.\label{eq:generaliseddifference}
 \end{equation}
  It may be called also an {$(\alpha,\beta)$-\emph{difference of order} $2s$}, or simply a  \emph{generalised difference}. In the case when $G={\mathbb T}$, we restrict $\alpha$ and $\beta$ to be in $\mathbb Z$ and ${\mathbb T}$ is identified with $[0,2\pi)$ . The vector space of functions in $ L^2(G)$ that can be expressed as some finite sum of $(\alpha,\beta)$-differences of order $2s$ is denoted by ${\mathcal D}_{\alpha,\beta,s}(G)$.  Thus, $f\in {\mathcal D}_{\alpha,\beta,s}({\mathbb R}) $ if and only if  there are $m\in {\mathbb N}$, $u_1,u_2,\ldots,u_m\in {\mathbb R}$  and $f_1,f_2,\ldots,f_m\in L^2({\mathbb R})$ such that
\[f=\sum_{j=1}^m\left[
\left(
e^{iu_j   \left(\frac{\alpha-\beta}{2}\right)  } +e^{-iu_j\left(\frac{\alpha-\beta}{2}\right)}
\right)
\delta_0-\left(
e^{iu_j   \left(\frac{\alpha+\beta}{2}\right)  } \delta_{u_j}+e^{-iu_j\left(\frac{\alpha+\beta}{2}\right)}\delta_{-u_j}
\right)
\right]^s\ast f_j.\]
  We prove  that if $f\in L^2({\mathbb R})$, $f\in {\mathcal D}_{\alpha,\beta,s}({\mathbb R})$ if and only if $\int_{-\infty}^{\infty}(x-\alpha)^{-2s}(x-\beta)^{-2s}|{\widehat f}(x)|^2<\infty$, in which case $f$ is a sum of $4s+1$ $(\alpha,\beta)$-differences of order $2s$. It follows that ${\mathcal D}_{\alpha,\beta,s}({\mathbb R})$ is the range of $(D^2-i(\alpha+\beta)-\alpha\beta I)^s({\mathbb R})$ on $W^{2s}({\mathbb R})$, and is a weighted $L^2$-space under the Fourier transform  in which the inner product  of $f$ and $g$ is $\int_{-\infty}^{\infty}\bigl(1+(x-\alpha)^{-2s}(x-\beta)^{-2s}\bigr){\widehat f}(x){\overline {{\widehat g}(x)}}\,dx$.
  
  Now,  if $\alpha, \beta\in {\mathbb Z}$, and if we take an  $(\alpha,\beta)$-difference $f$ in $L^2({\mathbb T})$ as in (\ref{eq:generaliseddifference}) with $G={\mathbb T}$, then    ${\widehat f}(\alpha)={\widehat g}(\beta)=0$. In \cite{nillsen3} it is proved that if $f$ in $L^2({\mathbb T})$ and ${\widehat f}(\alpha)={\widehat g}(\beta)=0$, then $f$ is a sum of $4s+1$ $(\alpha,\beta)$-differences of order $s$. Thus, the results obtained here extend the results obtained  in \cite{nillsen3}, for the compact case of the circle group ${\mathbb T}$, to the non-compact case of ${\mathbb R}$.   The techniques used here develop  the approach in \cite{nillsen3}, so as to deal with the additional complexities in going from $\mathbb T$ to  $\mathbb R$.

\section{Further notations and background}
\setcounter{equation}{0}
First we  need some  notions relating to partitions of an interval.  

{\bf Definitions.} If $J$ is an interval, $\lambda(J)$ denotes its length. A \emph{closed-interval partition} is a sequence $R_0,R_1,\ldots,R_{r-1}$ of closed  intervals of positive length such that $r=1$ or, when $r\ge 2$, the right hand endpoint of $R_j$ is the left hand endpoint of $R_{j+1}$ for all $j=0,1,\ldots,r-2$. We may refer to such  a closed-interval partition  as $\{R_0,R_1,\ldots,R_{r-1}\}$,  where we understand that the intervals $R_j$ may be rearranged so as to get a sequence forming a closed-interval partition. In this case if   we put  $J=\cup_{j=0}^{r-1}R_j$  then $\{R_0,R_1,\ldots,R_{r-1}\}$ may be called a \emph{closed-interval partition of J}.  If  $\{R_0,R_1,\ldots,R_{r-1}\}$ is a closed-interval partition  and $\{S_0,S_1,\ldots,S_{s-1}\}$ is also a  closed-interval partition, then   $\{R_j\cap S_k:0\le j\le r-1, 0\le k\le s-1\ {\rm and}\ \lambda(R_j\cap S_k)>0\}$ is  a closed-interval partition,  and we call it the \emph{refinement} of the closed-interval partitions $\{R_0,R_1,\ldots,R_{r-1}\}$ and $\{S_0,S_1,\ldots,S_{s-1}\}$. Finally, if $A$ is a set, $A^c$ will denote its complement.

 \begin{lemma}\label{lemma:partitionestimate} Let $J$ be a closed interval with $\lambda(J)>0$. Let $R_0,R_1,\ldots, R_{r-1}$ be  $r$  intervals in a closed-interval partition  of a closed interval $J$. Let $S_0,S_1,\ldots,S_{s-1}$  be $s$  intervals in a closed-interval partition  of a closed interval $K$, and assume that $\lambda(R_j\cap S_k)\ne\emptyset$ for at least one pair $j,k$.   Then the refinement of $R_0,R_1,\ldots, R_{r-1}$ and $S_0,S_1,\ldots,S_{s-1}$ is a closed-interval partition of $J\cap K$ and it has at most $r+s-1$ elements.
 \end{lemma}
 {\bf Proof.} It is easily checked that  the refinement of $R_0,R_1,\ldots, R_{r-1}$ and $S_0,S_1,\ldots,S_{s-1}$ is a closed-interval partition of $J\cap K$. Now, for any $r,s$ and  closed-interval partitions ${\mathcal P}=\{R_0,R_1,\ldots, R_{r-1}\}$ and ${\mathcal Q}=\{S_0,S_1,\ldots,S_{s-1}\}$ let's put
  \begin{align*}
 {\mathcal A}({\mathcal P},{\mathcal Q})&=\{(j, k):0\le j\le r-1, 0\le k\le s-1\ {\rm and}\ \lambda(R_j\cap S_k)>0\},\ {\rm and}\\
  {\mathcal B}({\mathcal P},{\mathcal Q})&=\{k:  0\le k\le s-1\ {\rm and}\ \lambda(R_r\cap S_k)>0\}. \\
      \end{align*}
      \vskip -0.7cm
      Let
      ${\overline{{\mathcal A}({\mathcal P},{\mathcal Q})}}$ and  ${\overline{{\mathcal B}({\mathcal P},{\mathcal Q})}}$ denote the number of elements in ${\overline{{\mathcal A}({\mathcal P}, {\mathcal Q})}}$ and ${\overline{{\mathcal B}({\mathcal P}, {\mathcal Q})}}$  respectively. The statement in the lemma is thus equivalent to saying that ${\overline{{\mathcal A}({\mathcal P},{\mathcal Q})}}\le r+s-1.$
  If $r=1$, we have $J=R_0$ and we see that ${\overline{{\mathcal A}({\mathcal P},{\mathcal Q})}}\le s=1+s-1$, so  in this case  the result holds  for any closed-interval partition $S_0,S_1,\ldots,S_{s-1}$.  We proceed by induction on $r$, by assuming that, for some  given $r\ge 2$,  for every  closed-interval partition ${\mathcal P}=\{R_0,R_1,\ldots, R_{r-1}\}$ and for any closed-interval partition ${\mathcal Q}$ having $s$ elements for arbitrary $s\in {\mathbb N}$,  we have 
  ${\overline{{\mathcal A}({\mathcal P},{\mathcal Q})}}\le r+s-1.$

       Now consider  closed-interval partitions ${\mathcal P}^{\prime}=\{R_0,R_1,\ldots,R_r\}$   and ${\mathcal Q}^{\prime}=\{S_0,S_1,\ldots,S_{s-1}\}$. We may assume that  ${\overline {{\mathcal A}({\{R_0, R_1,\ldots,R_{r-1}\}, {\mathcal Q}^{\prime}) }}}\ge 1$ , for otherwise we have \hfill\break ${\overline {{\mathcal A}({\{R_0, R_1,\ldots,R_{r-1}\}, {\mathcal Q}^{\prime}) }}}=0 $, and then
       \[{\overline {{\mathcal A}({\{R_0, R_1,\ldots,R_{r}\}, {\mathcal Q}^{\prime}) }}}={\overline {{\mathcal A}({\{R_{r}\}, {\mathcal Q}) }}}\le s\le (r+1)+s-1,\]
        as we have seen the lemma is true when one of the partitions has a single element. That is, in this case, the truth of the lemma for $r$ implies the truth of the lemma for $r+1$.
       
        Now, when ${\overline {({\mathcal A}({\{R_0, R_1,\ldots,R_{r-1}\}, {\mathcal Q}^{\prime}) }}}\ge 1$, let $s_1\in \{1,2,\ldots,s-1\}$ be the  maximum of all $k\in    \{1,2,\ldots,s-1\}$ such that $\lambda(A_j\cap B_k)>0$ for some $j\in \{0,1,2\ldots,r-1\}$.
 By the inductive assumption, ${\overline {{\mathcal A}{(\{R_0, R_1,\ldots,R_{r-1}\}, {\mathcal Q}^{\prime}) }}}\le r+s_1-1$.
  Also, as we go from $r$ to $r+1$, the single interval $R_r$ is adjoined to $R_0,R_1,\ldots,R_{r-1}$ on the right. So, we see that $\overline {
  {\mathcal B}({\mathcal P}^{\prime},{\mathcal Q}^{\prime})}\le s-(s_1-1)$. As well, it is clear that 
  \[{\overline { {\mathcal A}({\mathcal P}^{\prime},{\mathcal Q}^{\prime})}}
  \le{\overline{{\mathcal A}(\{R_0,R_1,\ldots,R_{r-1}\}, {\mathcal Q}^{\prime})}}+\overline {
  {\mathcal B}({\mathcal P}^{\prime},{\mathcal Q}^{\prime})}  .\]
    Using the inductive assumption it follows that 
     \[{\overline { {\mathcal A}({\mathcal P}^{\prime},{\mathcal Q}^{\prime})}}
  \le r+s_1-1+s-(s_1-1)=(r+1)+s-1,\]
 and we see that  assuming  the inductive assumption holds  for $r$   implies that it holds for $r+1$. Invoking induction completes the proof.\hfill$\square$

 \begin{lemma}\label{lemma:quadratics}

 Let $a,b,c,d\in {\mathbb R}$ with $c<d$,  $a\le  d$ and  $b\ge c$.   Put
 
 (i) $a^{\prime}=a$ if $a\in [c,d]$, and  $a^{\prime}=c$ if $a<c$; 
 
 (ii) $b^{\prime}=b$ if $b\in [c,d]$, and  $b^{\prime}=d$ if $d<b$.  
 
 Then, 
 \[|(u-a)(u-b)|\ge |(u-a^{\prime})(u-b^{\prime})|,\ \hbox{for all}\  u\in [c,d].\]
 \end{lemma}
{\bf Proof.} If $a^{\prime}=a$ or $b^{\prime}=b$, the result is easily checked. The only other case is when $a^{\prime}=c$ and $b^{\prime}=d$. In this case we have
\[(u-a)(b-u)-(u-c)(d-u)=(a+b-c-d)u-ab+cd,\]
and this expression is linear in $u$ and non-negative for $u=c$ and $u=d$. Consequently,
$(u-a)(b-u)\ge (u-c)(d-u)$ for all $ u\in [c,d]$, and the proof is complete.
\hfill$\square$

 \begin{lemma} \label{lemma:multi}
  Let $s,m\in {\mathbb N}$ with $m\ge 4s+1$, and   let $V_1,V_2,\ldots,V_m$ be closed and bounded subintervals of $\mathbb R$. Let $c_1,c_2,\ldots,c_m, d_1,d_2,\ldots,d_m\in {\mathbb R}$ be such that $c_j,d_j\in V_j$ for all $j=
  1,2,\ldots,m$. Then, there is a number $M>0$, depending upon $s$ and $m$ only and independent of the $V_j, c_j$ and $d_j$, such that
 \[ \int_{\textstyle{\prod_{t=1}^mV_t }}\, \,\frac{du_1du_2\ldots du_m}{\displaystyle \sum_{t=1}^{m}\left(u_t-c_t\right) ^{2s}\left(u_t-d_t \right)^{2s}}
\ \le\  M\,\Big(\max\Big\{\mu(V_1),\mu(V_2),\ldots,\mu(V_m)\Big\}  \Big)^{m-4s}.\]
  \end{lemma}
{\bf Proof.} See \cite[Lemma 4.2]{nillsen3}.\hfill$\square$

 \begin{theorem} \label{theorem:characterisation} Let $f\in L^2({\mathbb R})$ and let $\mu_1,\mu_2,\ldots,\mu_r\in M({\mathbb R})$. Then the following conditions  (i) and (ii) are equivalent.
 
 (i) There are   $f_1,f_2,\ldots, f_r\in L^2({\mathbb R})$ such that 
 $f=\sum_{j=1}^r\mu_j\ast f_j.$
 \vskip0.2cm
(ii) \hskip 3.6cm$\displaystyle \int_{-\infty}^{\infty}\,\frac{|{\widehat f }(x)|^2}{\displaystyle{\sum_{j=1}^r}|{\widehat \mu_j}(x)|^2}\,dx\ <\ \infty.$
\end{theorem}

{\bf Proof.} This is essentially proved in \cite[pages 411-412 ]{meisters1}, but see also \cite[pages 77-88]{nillsen1} and \cite[page 23]{nillsen2}. \hfill$\square$

\section{Main results}
\setcounter{equation}{0}
The main aim in this paper is to prove the following result.
\begin{theorem} \label{theorem:main} Let $s\in {\mathbb N}$ and let  $\alpha,\beta\in {\mathbb R}$. Let ${\mathcal D}_{\alpha,\beta,s}({\mathbb R})$ be the  vector space  of functions $f\in L^2({\mathbb R})$ that can be expressed as some finite sum of $(\alpha,\beta)$-differences of order $2s$.   Then the following conditions (i) - (iii) are equivalent for a function $f\in L^2({\mathbb R})$.
\vskip 0.2cm
(i) $\displaystyle \int_{-\infty}^{\infty}\frac{|{\widehat f}(x)|^2}{(x-\alpha)^{2s}(x-\beta)^{2s}}dx<\infty.$
\vskip 0.2cm
(ii) $f\in {\mathcal D}_{\alpha,\beta,s}({\mathbb R})$.

(iii) There are   $u_1,u_2,\ldots,u_{4s+1}\in {\mathbb R}$  and $f_1,f_2,\ldots,f_{4s+1}\in L^2({\mathbb R})$ such that
\begin{equation}f=\sum_{j=1}^{4s+1}\left[
\left(
e^{iu_j   \left(\frac{\alpha-\beta}{2}\right)  } +e^{-iu_j\left(\frac{\alpha-\beta}{2}\right)}
\right)
\delta_0-\left(
e^{iu_j   \left(\frac{\alpha+\beta}{2}\right)  } \delta_{u_j}+e^{-iu_j\left(\frac{\alpha+\beta}{2}\right)}\delta_{-u_j}
\right)
\right]^s\ast f_j.\label{eq:expressed}
\end{equation}

When the preceding conditions hold for a given function $f\in L^2({\mathbb R})$, for almost all\hfill\break  $(u_1,u_2,\ldots,u_{4s+1})\in {\mathbb R}^{4s+1}$, there are $f_1,f_2,\ldots,f_{4s+1}\in L^2({\mathbb R})$ such that (\ref{eq:expressed}) holds. Also, ${\mathcal D}_{\alpha,\beta,s}({\mathbb R})$   is a Hilbert space with the inner product ${\langle\,,\rangle}_{\alpha,\beta,s}$ given by
\[{\langle f,g\rangle}_{\alpha,\beta,s} =\int_{-\infty}^{\infty}\left(1+\frac{1}{(x-\alpha)^{2s}(x-\beta)^{2s}}\right)  {\widehat f}(x){\overline{{\widehat g}(x)}}dx,\ {\rm for}\ f,g\in {\mathcal D}_{\alpha,\beta,s}({\mathbb R}).\]
The operator $(D^2-(\alpha+\beta)D-\alpha\beta I)^s$ is  a linear, bounded and invertible operator that maps $W^{2s}({\mathbb R})$
onto ${\mathcal D}_{\alpha,\beta,s}({\mathbb R})$.
\end{theorem}

{\bf Proof.} If (iii) holds it is clear that (ii) also holds.

Let (ii) hold. If $u\in {\mathbb R}$, define $\lambda_u\in M({\mathbb R})$ by
\begin{equation}\lambda_u=\frac{1}{2}\left[e^{iu\left(\frac{\alpha-\beta}{2}\right)}+e^{-iu\left(\frac{\alpha-\beta}{2}\right) } \right]\delta_0- \frac{1}{2}\left[\  e^{iu\left(\frac{\alpha+\beta}{2}\right)}\,\delta_{u}+e^{-iu\left(\frac{\alpha+\beta}{2}\right)}\,\delta_{-u}\right].\label{eq:lambdasubb}
\end{equation}
The Fourier transform ${\widehat {\lambda}_u}$ of $\lambda_u$ is given for $x\in {\mathbb R}$ by
\begin{equation}
{\widehat {\lambda}}_u(x)=\sin\left(u(x-\alpha)\right) \sin\left(u(x-\beta)\right). \label{eq:Fouriertransform}
\end{equation}
So, if $u\in {\mathbb R}$ and  $f, g\in L^2({\mathbb R})$ are such that $f=\lambda_u^s\ast g$, we have
\begin{align*}\int_{-\infty}^{\infty}\frac{|{\widehat f}(x)|^2}{(x-\alpha)^{2s}(x-\beta)^{2s}}dx&=\int_{-\infty}^{\infty}\frac{\sin^{2s}\left(u(x-\alpha)\right) \sin^{2s}\left(u(x-\beta)\right)}{(x-\alpha)^{2s}(x-\beta)^{2s}}|{\widehat g}(x)|^2dx<\infty. \end{align*}
Using the definition of $\lambda_u$ in (\ref{eq:lambdasubb}), we    deduce that (ii) implies (i).

Now, we assume that (i) holds, and we will prove that (iii) holds.   Let  $c>0$ be given and let  $x\in \{\mathbb R\}$ also be given   but with $x\notin \{\alpha,\beta\}$.   Put, for each $k\in {\mathbb Z}$,
\begin{equation}a_k=\frac{k\pi}{|x-\alpha|}, b_k=\frac{k\pi}{|x-\beta|}, a_k^{\prime}=\frac{(k-1/2)\pi}{|x-\alpha|}  \ {\rm and}\  b_k^{\prime}=\frac{(k-1/2)\pi}{|x-\beta|}. \label{eq:zeros}
\end{equation}
Then put, again for each $k\in {\mathbb Z}$,
\begin{equation}A_k=[a_k^{\prime},a_{k+1}^{\prime}]\ {\rm and }\ B_k=[b_k^{\prime},b_{k+1}^{\prime}].\label{eq:intervals}
\end{equation}
Note that $a_k$ is the mid-point of $A_k$ and $b_k$ is the mid-point of $B_k$. The points $a_k$ are the zeros of $u\longmapsto \sin (u(x-\alpha))$, while the $b_k$ are the zeros of $u\longmapsto \sin (u(x-\beta))$.
It is immediate from the definitions that, for each $k\in {\mathbb Z}$,
\begin{equation}\lambda(A_k)= \frac{\pi}{|x-\alpha|}\ {\rm and}\ \lambda(B_k)= \frac{\pi}{|x-\beta|}.\label{eq:intervalslength}
\end{equation}
We will use the notation that $d_{\mathbb Z}(x)$ denotes the distance from $x\in {\mathbb R}$ to the nearest integer. Note that $d_{\mathbb Z}(x)=|x|$ if and only if $-1/2\le x\le1/2$. Note also that $|\sin(\pi x)|\ge 2d_{\mathbb Z}(x)$ for all $x\in {\mathbb R}$ (for example see \cite[page 89]{nillsen2} or \cite[page 233]{stromberg1}).

Now 
\begin{align*}u\in A_j&\Longrightarrow
\frac{(j-1/2)\pi}{|x-\alpha|}\le u\le\frac{(j+1/2)\pi}{|x-\alpha|}
\Longrightarrow -1/2\le |x-\alpha|\left| \frac{u}{\pi}-\frac{j}{|x-\alpha|}\right|\le1/2.
\end{align*}
So, for $u\in A_j$,
\begin{align}
|\sin(u(x-\alpha))|
&=\left|\sin\left( \pi|x-\alpha|\left|\frac{u}{\pi}-\frac{j}{|x-\alpha|}\right|\right) \right|\nonumber\\
&\ge 2d_{\mathbb Z}\left( |x-\alpha|\left|\frac{u}{\pi}-\frac{j}{|x-\alpha|}\right|\right)\nonumber\\
&= 2|x-\alpha|\left|\frac{u}{\pi}-\frac{j}{|x-\alpha|}\right|\nonumber\\
&=\frac{2}{\pi}|x-\alpha|\left|u-\frac{j\pi}{|x-\alpha}\right|.\label{eq:sine1}
\end{align}
Similarly, for $u\in B_k$,
\begin{equation}|\sin(u(x-\beta))|\ge\frac{2}{\pi}|x-\beta|\left|u-\frac{k\pi}{|x-\beta|}\right|.\label{eq:sine2}
\end{equation}
So, for $u\in A_j\cap B_k$ we have
\[|\sin(u(x-\alpha))\sin(u(x-\beta))|\ge\frac{4}{\pi^2}|(x-\alpha)(x-\beta)|\left|u-\frac{j\pi}{|x-\alpha}\right|\cdot\left|u-\frac{k\pi}{|x-\beta}\right|.\]
That is, for $u\in  A_j\cap B_k$ we have
\begin{equation}
|\sin(u(x-\alpha))\sin(u(x-\beta))|\ge \frac{4}{\pi^2}|(x-\alpha)(x-\beta)|\cdot |u-a_j|\cdot |u-b_k|,\label{eq:sineinequality}
\end{equation}
where $a_j$ and $b_k$ are the points as given in (\ref{eq:zeros}), with $a_j$ the midpoint of $A_j$ and $b_k$ the midpoint of $B_k$.

 Now let the intervals $A_j$ such that $\lambda(A_j\cap [-c,c]\,)>0$ be 
$ A_{m_1},\ldots, A_{m_1+r-1}$,
and  let the intervals $B_k$ such that $\lambda(B_k\cap [-c,c]\,)>0$ be
$  B_{m_2},\ldots, B_{m_2+s-1}.$
Put
\begin{align}A_j^{\prime}&=A_{m_1+j},\  \hbox{for $j=0,1,\ldots, r-1$},\ \hbox{and put}\nonumber\\
B_k^{\prime}&=B_{m_2+k},\  \hbox{for $k=0,1,\ldots, s-1$}.\label{eq:ABprime}
\end{align}
Then, putting
\begin{equation}
{\mathcal P}_1=\{A_0^{\prime}, A_1^{\prime},\ldots,A_{r-1}^{\prime}\} \ {\rm and}\ {\mathcal P}_2=\{B_0^{\prime}, B_1^{\prime},\ldots,B_{s-1}^{\prime}\},\label{eq:partitions}
\end{equation}
we see that ${\mathcal P}_1$ and ${\mathcal P}_2$ are closed-interval partitions.      If we put 
 \begin{equation}
 {\mathcal A}=\{(j,k):0\le j\le r-1, 0\le k\le s-1,\ {\rm and}\  \lambda(A_j^{\prime}\cap B_k^{\prime})>0\},\nonumber
 \end{equation}
 and if we let ${\mathcal P}$ be the refinement of ${\mathcal P}_1$ and ${\mathcal P}_2$, 
 we have
 \begin{equation}
 {\mathcal P}=\{A_j^{\prime}\cap B_k^{\prime}: (j,k)\in {\mathcal A}\}\ {\rm and}\ [-c,c]\subseteq \bigcup_{(j,k)\in {\mathcal A}}A_j^{\prime}\cap B_k^{\prime}. \label{eq:refinement}
 \end{equation}
 
 Now, from (\ref{eq:intervalslength}) we see that all lengths of the  $r$ intervals in the closed-interval partition ${\mathcal P}_1$  equal $\pi/|x-\alpha|$, so that $(r-2)\pi/|x-\alpha|< 2c.$
Hence,
  \begin{equation}  1\le r< \frac{2c|x-\alpha|}{\pi}+2=\frac{2c}{\pi}\left(1+\frac{\pi}{c|x-\alpha|}\right)|x-\alpha|. \label{eq:inequalitya}
  \end{equation}
    
 Let $0<\delta<1/2$. Then, if  $|x-\alpha|>\pi\delta/c$, we have from (\ref{eq:inequalitya}) that
 \begin{equation}1\le r< \frac{2c}{\pi}\left(1+\frac{1}{\delta }\right)|x-\alpha|.\label{eq:alphahigh}
 \end{equation}
 On the other hand, if  $|x-\alpha|\le\pi\delta/c$, as $0<\delta<1/2$ we have $2c<\pi/|x-\alpha|$, and it follows from (\ref{eq:intervalslength}) that $[-c,c]\subseteq A_0$, so that  $m_0 =0$ and  
 \begin{equation}r=1. \label{eq:alphalow}
 \end{equation}

Again let  $0<\delta<1/2$. Then, as in the preceding argument, but with $\beta$ replacing $\alpha$,  if   $|x-\beta|>\pi\delta/c$ we have
   \begin{equation}1\le s< \frac{2c}{\pi}\left(1+\frac{1}{\delta }\right)|x-\beta|,\label{eq:betahigh}
    \end{equation} 
 while if  $|x-\beta|\le\pi\delta/c$, we have
   \begin{equation}
   s=1.\label{eq:betalow}
    \end{equation}
  
Now, again let  $0<\delta<1/2$. Assume that either $|x-\alpha|>\pi\delta/c$ or    $|x-\beta|>\pi\delta/c$, with both perhaps holding. In the case that   $|x-\alpha|>\pi\delta/c$ and $|x-\beta|\le\pi\delta/c$, we have from (\ref{eq:alphahigh}) and   (\ref{eq:betalow}) 
  that
\begin{align}r+s-1\le2\max\{r,s\}
&<2\max\left\{\frac{2c}{\pi}\left(1+\frac{1}{\delta }\right)|x-\alpha|, 1\right\}\nonumber\\
&\le 2\max\left\{\frac{2c}{\pi}\left(1+\frac{1}{\delta }\right)|x-\alpha|, \frac{2c}{\pi}\left(1+\frac{1}{\delta }\right)|x-\beta|\right\}\nonumber\\
&=\frac{4c}{\pi}\left(1+\frac{1}{\delta }\right)\max\{|x-\alpha|,|x-\beta|)\}.\label{eq:estimateb}
\end{align}
The argument  that produced (\ref{eq:estimateb}) is symmetric in $\alpha$ and $\beta$, and we see from (\ref{eq:alphahigh}), (\ref{eq:alphalow}), (\ref{eq:betahigh}) and (\ref{eq:betalow}) that in all cases when either $|x-\alpha|>\pi\delta/c$ or    $|x-\beta|>\pi\delta/c$ we have
 \begin{equation}r+s-1\le  2\max\{r,s\}\le \frac{4c}{\pi}\left(1+\frac{1}{\delta }\right)\max\Big\{|x-\alpha|,|x-\beta|\Bigr\}.\label{eq:partitionestimate1}
 \end{equation}
 
Also, observe that if $0<\delta<1/2$, $|x-\alpha|\le\pi\delta/c$ and    $|x-\beta| \le\pi\delta/c$, we have from (\ref{eq:alphalow}) and (\ref{eq:betalow}) that
 \begin{equation}
r=s=1.\label{eq:partitionestimate2}
 \end{equation}
Note that in the above,  $a_k,b_k,A_k, B_k$ and so on, depend upon $x$. In particular, $r$ and $s$ depend upon $x$.

 We now take $m\in {\mathbb N}$ with $m\ge 4s+1$, and we estimate the integral
\[\int_{[-c,c]^m}\frac{du_1du_2\ldots du_m}{\displaystyle \sum_{j=1}^m\sin^{2s}u_j\left(x-\alpha \right)\, \sin^{2s}u_j\left( x-\beta\right)},\] 
allowing for the  different values $x$ may be, but recall that $x\notin \{\alpha,\beta\}$. We let ${\mathcal P}_1$, ${\mathcal P}_2$ be the closed-interval partition   as given in (\ref{eq:partitions}) and let  ${\mathcal P}$ be their refinement as given in (\ref{eq:refinement}). 
We have, using  the definitions and  (\ref{eq:zeros}), (\ref{eq:sineinequality}), (\ref{eq:ABprime}) and (\ref{eq:refinement}),
\begin{align}
&\int_{[-c,c]^m}\frac{du_1du_2\ldots du_m}{\displaystyle \sum_{j=1}^m\sin^{2s}u_j\left(x-\alpha \right)\, \sin^{2s}u_j\left( x-\beta\right)}\nonumber\\
&\le \sum_{(j_1,k_1), (j_2,k_2),\ldots,(j_m,k_m)\in {\mathcal A}} \int_{\textstyle\prod_{t=1}^m A_{j_t}^{\prime} \cap B_{k_t}^{\prime}}\frac{du_1du_2\ldots du_m}{\displaystyle \sum_{j=1}^m\sin^{2s}u_j\left(x-\alpha \right)\, \sin^{2s}u_j\left( x-\beta\right)}\nonumber\\
&\le{ \frac{\pi^{4s}}{2^{4s}(x-\alpha)^{2s}(x-\beta)^{2s}} }\sum_{(j_1,k_1),  \ldots,(j_m,k_m)\in {\mathcal A}}\ \int_{\textstyle\prod_{t=1}^m A_{j_t}^{\prime} \cap B_{k_t}^{\prime}}\frac{du_1du_2\ldots du_m}{ (u_j-a_{m_1+j_t})^{2s} (u_j-b_{m_2+k_t})^{2s}}. \label{eq:fundamental1}
\end{align}
Now in (\ref{eq:fundamental1}), the points $a_{m_1+j_t},b_{m_2+k_t}$ do not necessarily belong to $A_{j_t}^{\prime} \cap B_{k_t}^{\prime}$. However, suppose that $a_{m_1+j_t},b_{m_2+k_t}\notin A_{j_t}^{\prime} \cap B_{k_t}^{\prime}$ with   $a_{m_1+j_t}\le b_{m_2+k_t}$ and with both $a_{m_1+j_t}$ and $b_{m_2+k_t}$  lying to the left of $A_{j_t}^{\prime} \cap B_{k_t}^{\prime}$. Let $y$ be the left endpoint of $A_{j_t}^{\prime} \cap B_{k_t}^{\prime}$.
Then,  $a_{m_1+j_t}\in A_{m_1+j_t}$ so we see that $[a_{m_1+j_t}, y]\subseteq A_{m_1+j_t}$. Similarly,  $[b_{m_2+k_t}, y]\subseteq B_{m_2+k_t}$. We deduce that 
\[[b_{m_2+k_t}, y]\subseteq  A_{m_1+j_t}\cap  B_{m_2+k_t}= A_{j_t}^{\prime}\cap B_{k_t}^{\prime},\]
so that $b_{m_2+k_t}\in A_{j_t}^{\prime}\cap B_{k_t}^{\prime}$,  and this contradicts the fact that  $b_{m_2+k_t}\notin A_{j_t}^{\prime} \cap B_{k_t}^{\prime}$.
 This argument, repeated for other cases, means we can say that if $a_{m_1+j_t},b_{m_2+k_t}\notin A_{j_t}^{\prime} \cap B_{k_t}^{\prime}$, then $a_{m_1+j_t}$ is on the left of $A_{j_t}^{\prime} \cap B_{k_t}^{\prime}$ and $b_{m_2+k_t}$ is on the right of $A_{j_t}^{\prime} \cap B_{k_t}^{\prime}$, or \emph{vice versa}. 
 
   Now,  let $t\in\{1,2,\ldots,m\}$ and  $(j_t,k_t)\in {\mathcal A}$ be given. In  Lemma \ref{lemma:quadratics}, let's take   $[c,d]$ to be $A_{j_t}^{\prime} \cap B_{k_t}^{\prime}$,  $a$ to be $\min\{a_{m_1+j_t}, b_{m_2+k_t}\}$, and $b$ to be $\max\{a_{m_1+j_t}, b_{m_2+k_t}\}$.  The observation at the end of the preceding paragraph implies that  we must have $a\le d$ and $b\ge c$. Also, note that either $a=a_{m_1+j_t}$ and $b=b_{m_2+k_t}$ or \emph{vice versa}. Then, from Lemma \ref{lemma:quadratics} we have
\begin{equation}
|(u-a_{m_1+j_t})(u-b_{m_2+k_t})|\ge |(u-a_{m_1+j_t}^{\prime})(u-b_{m_2+k_t}^{\prime})|,\ \hbox{for all}\ u\in A_{j_t}^{\prime} \cap B_{k_t}^{\prime}. \label{eq:fundamental2}
\end{equation}

  Now let $0<\delta<1/2$ and assume that we have either $|x-\alpha|>\pi\delta/c$ or $|x-\beta|>\pi\delta/c$. Then from Lemma \ref{lemma:partitionestimate}, the  right hand side of  (\ref{eq:partitionestimate1}) gives an upper bound  for the number of elements in $\mathcal P$.  Using (\ref{eq:partitionestimate1}), and using   (\ref{eq:sineinequality}), (\ref{eq:refinement}), (\ref{eq:partitionestimate2}), (\ref{eq:fundamental1}), (\ref{eq:fundamental2}),  Lemma \ref{lemma:multi} and the assumption that $m\ge 4s+1$,  we have in this case that
   \begin{align}
   &\int_{[-c,c]^m}\frac{du_1du_2\ldots du_m}{\displaystyle \sum_{j=1}^m\sin^{2s}u_j\left(x-\alpha \right)\, \sin^{2s}u_j\left( x-\beta\right)}\nonumber\\
&\le{ \frac{\pi^{4s}}{2^{4s}(x-\alpha)^{2s}(x-\beta)^{2s}} }\sum_{(j_1,k_1), \ldots,(j_m,k_m)\in {\mathcal A}}\ \int_{\textstyle\prod_{t=1}^m A_{j_t}^{\prime} \cap B_{k_t}^{\prime}}\frac{du_1du_2\ldots du_m}{ |u_j-a_{m_1+j_t}^{\prime}|^{2s} |u_j-b_{m_2+k_t}^{\prime}|^{2s}}\nonumber
\end{align}
 
\begin{align}
&\le{ \frac{\pi^{4s}M}{2^{4s}(x-\alpha)^{2s}(x-\beta)^{2s}} }\sum_{(j_1,k_1),  \ldots,(j_m,k_m)\in {\mathcal A}}\,\Big(\max\Big\{\lambda(A_{j_1}^{\prime} \cap B_{k_1}^{\prime}),\ldots,\lambda(A_{j_m}^{\prime} \cap B_{k_m}^{\prime}) \Big\}\Big)^{m-4s},\nonumber\\
& \hskip 8.6cm\hbox{where $M$ is the constant in  Lemma\   \ref{lemma:multi} },\nonumber\\
&\le{ \frac{\pi^{4s-m}2^{2m-4s}c^m(\delta+1)^mM}{\delta^m(x-\alpha)^{2s}(x-\beta)^{2s}} }\max\Big\{|x-\alpha|^m,|x-\beta|^m\Bigr\}\min\left\{ \frac{\pi^{m-4s}}{|x-\alpha|^{m-4s}},  \frac{\pi^{m-4s}}{|x-\beta|^{m-4s}}\right\},\nonumber\\
& \hskip 13.6cm\hbox{using   (\ref{eq:intervalslength}),}\nonumber\\
&\le Q\max\left\{\, \frac{(x-\alpha)^{2s}}{(x-\beta)^{2s}}\,,\,\frac{(x-\beta)^{2s}}{(x-\alpha)^{2s}}\,\right\},\label{eq:integralestimate1}
 \end{align}
 for all $x\notin\{\alpha,\beta\}$ with either $|x-\alpha|>\pi\delta/c$ or $|x-\beta|>\pi\delta/c$. Note that the constant $Q$ in (\ref{eq:integralestimate1}) is independent of $x$.

 CASE I: $\alpha\ne \beta$.
 
 In this case, choose $\delta$ so that
\begin{equation}0<\delta<\min\left\{\frac{1}{2},\frac{c|\alpha-\beta|}{2\pi}\right\}.\label{eq:delta1}
\end{equation}
Then, define disjoint intervals $J, K$ by putting
 \[J=\left[\alpha-\frac{\pi \delta}{c}, \alpha+\frac{\pi\delta}{c}\right]\ {\rm and}\ K=\left[\beta-\frac{\pi \delta}{c}, \beta+\frac{\pi\delta}{c}\right]. \]
 Clearly, there is $C_1>0$ such   that 
 
\begin{equation}
\max\left\{\, \frac{(x-\alpha)^{2s}}{(x-\beta)^{2s}}\,,\,\frac{(x-\beta)^{2s}}{(x-\alpha)^{2s}}\,\right\}\le C_1,\ \hbox{for all}\ x\in (J\cup K)^c.\label{eq:X}
\end{equation}
 As well, $(x-\beta)^{-2s} $ is bounded on $J$, so we see that there is $C_2>0$ such that 
 \begin{equation}
\max\left\{\, \frac{(x-\alpha)^{2s}}{(x-\beta)^{2s}}\,,\,\frac{(x-\beta)^{2s}}{(x-\alpha)^{2s}}\,\right\}(x-\alpha)^{2s}\le C_2,\ \hbox{for all}\ x\in J\cap\{\alpha\}^c.\label{eq:Y}
\end{equation}
And, as 
 $(x-\alpha)^{-2s} $ is bounded on $K$,   there is $C_3>0$ such that
  \begin{equation}
\max\left\{\, \frac{(x-\alpha)^{2s}}{(x-\beta)^{2s}}\,,\,\frac{(x-\beta)^{2s}}{(x-\alpha)^{2s}}\,\right\}(x-\beta)^{2s}\le C_3,\ \hbox{for all}\ x\in K\cap\{\beta\}^c.\label{eq:Z}
\end{equation}
 
 We now  have from (\ref{eq:integralestimate1}), (\ref{eq:X}),  (\ref{eq:Y}) and  (\ref{eq:Z}), that  
 \begin{align}
& \int_{-\infty}^{\infty}\left(\int_{[-c,c]^m}\frac{du_1du_2\ldots du_m}{\displaystyle \sum_{j=1}^m\sin^{2s}u_j\left(x-\alpha \right)\, \sin^{2s}u_j\left( x-\beta\right)}\right)\,|{\widehat f}(x)|^2\,dx\nonumber\\
&\le C_1Q\int_{(J\cup K)^c}|{\widehat f}(x)|^2\,dx+C_2Q\int_J\frac{|{\widehat f}(x)|^2}{(x-\alpha)^{2s}}dx+C_3Q\int_K\frac{|{\widehat f}(x)|^2}{(x-\beta)^{2s}}dx\nonumber\\
&<\infty,\label{eq:conclusion1}
\end{align}
as  we are assuming that $\int_{-\infty}^{\infty}|{\widehat f}(x)|^2(x-\alpha)^{-2s}(x-\beta)^{-2s}<\infty$.  

CASE II. $\alpha=\beta$.   

Let's assume that $\alpha\in (-c, c)$ and that
\begin{equation}\delta<\min\left\{\frac{1}{2}, \frac{c(c-|\alpha|)}{\pi}\right\}. \label{eq:deltachoice}
\end{equation}
 Put $L=  ( \,\alpha-\pi\delta/c, \alpha+\pi\delta/c\,)$, and observe that because of (\ref{eq:deltachoice}), $L\subseteq (-c,c)$. 
 If $x\in L$, we have $|x-\alpha|<\pi\delta/c$ and it follows from (\ref{eq:refinement})  and (\ref{eq:alphalow}) that $r =1$. Now as $r=1$ and as $A_0\cap [-c,c]\ne \emptyset$,   we see that $[-c,c]\subseteq A_0=A_0^{\prime}= B_0^{\prime}$. In this case, in  (\ref{eq:sine1})  we must have $k=0$, and as $L\subseteq [-c,c]\subseteq A_0^{\prime}$,  we deduce  from (\ref{eq:sine1})  that
 \begin{equation}\ |\sin(u(x-\alpha))|\ge \frac{2}{\pi}|u|\cdot|x-\alpha|,\ \hbox{for all}\ x\in L\ \hbox{and}\ u\in (-c,c).\label{eq:refinementone}
 \end{equation}
 Let $C>0$ be such that 
 \begin{equation}
 \sum_{j=1}^mu_j^{4s}\ge C\left(\sum_{j=1}^mu_j^2\right)^{2s},\ \hbox{for all}\ (u_1,u_2,\dots,u_m)\in {\mathbb R}^m.\label{eq:euclideanconstant}
 \end{equation}
 We now  have from (\ref{eq:refinementone}) and (\ref{eq:euclideanconstant}) that if $m\ge 4s+1$ and  $x\in L$,
\begin{align}
&\int_{[-c,c]^m}\frac{du_1du_2\ldots du_m}{\displaystyle \sum_{j=1}^m\sin^{4s}u_j\left(x-\alpha \right)  }\nonumber\\
  &\le\frac{\pi^{4s}}{2^{4s}(x-\alpha)^{4s}}\ \int_{[-c,c]^m} \frac{du_1du_2\ldots du_m}{\displaystyle\sum_{j=1}^mu_j^{4s}}\nonumber
  \end{align}

  \begin{align}
        &\hskip 1.2cm\le \frac{1}{C}\cdot \frac{\pi^{4s}}{2^{4s}(x-\alpha)^{4s} }\int_{{[-c,c]^m}}\frac{du_1du_2\ldots du_m}{\displaystyle\left(\sum_{j=1}^mu_j^2\right)^{2s}}\nonumber\\
    &\hskip 1.2cm\le \frac{D}{C}\cdot\frac{\pi^{4s}}{2^{4s}(x-\alpha)^{4s} }\int_0^{c{\sqrt m}}r^{m-4s-1}\,dr,\nonumber\\
    &\hskip 7cm\hbox{for some $D>0$, by\ \cite[pages 394-395]{stromberg1}},\nonumber\\
    &\hskip 1.2cm\le \frac{G}{ (x-\alpha)^{4s}},\label{eq:GinequalityA}
\end{align}
for some $G>0$ that is independent of $x\in L\cap\{\alpha\}^c$. 

On the other hand,  if $x\notin L$ we have $|x-\alpha|\ge\pi\delta/c$, so that if we apply (\ref{eq:integralestimate1}) with $\alpha=\beta$ we have
\begin{equation}
\int_{[-c,c]^m}\frac{du_1du_2\ldots du_m}{\displaystyle \sum_{j=1}^m\sin^{4s}u_j\left(x-\alpha \right)}\,\le Q<\infty.\label{eq:GinequalityB}
\end{equation}
 Assuming that $|\alpha|<c$, we now have, using (\ref{eq:GinequalityA}) and (\ref{eq:GinequalityB}) and the fact that ${\mathbb R}=L\cup L^c$,
 \begin{align}
&\int_{-\infty}^{\infty}\left(\int_{[-c,c]^m}\frac{du_1du_2\ldots du_m}{\displaystyle \sum_{j=1}^m\sin^{4s}u_j\left(x-\alpha \right)\,}\right)|{\widehat f}(x)|^2 dx\nonumber\\
&\le G\int_{L}\frac{|{\widehat f}(x)|^2}{(x-\alpha)^{4s}} dx\ +Q \int_{L^c}|{\widehat f}(x)|^2 dx\ \ \nonumber\\
&<\infty,\label{eq:conclusion2}
\end{align}
as $\alpha=\beta$ and we are assuming that $\int_{-\infty}^{\infty}|{\widehat f}(x)|^2(x-\alpha)^{-2s}(x-\beta)^{-2s}\,dx<\infty$.

We have considered the cases $\alpha\ne \beta$ and $\alpha=\beta$. The d\'enoument results  from using Fubini's Theorem,  (\ref{eq:conclusion1}) and (\ref{eq:conclusion2}). We see that  provided  $|\alpha|<c$ and  $m\ge 4s+1$, in both cases  we have
\[\int_{[-c,c]^m}\left(\int_{-\infty}^{\infty}\frac{|{\widehat f}(x)|^2 dx}{\displaystyle \sum_{j=1}^m\sin^{2s}u_j\left(x-\alpha \right)\, \sin^{2s}u_j(x-\beta)}\right)du_1du_2\ldots du_m<\infty.\]
We conclude from this that for almost all $(u_1,u_2,\ldots,u_m)\in [-c,c]^m$, 
\begin{equation}\int_{-\infty}^{\infty}\frac{|{\widehat f}(x)|^2 dx}{\displaystyle \sum_{j=1}^m\sin^{2s}(u_j(x-\alpha)) \, \sin^{2s}(u_j(x-\beta))}
<\infty.\label{eq:almosteverywhere}
\end{equation}
By letting $c$ tend to $\infty$ through a sequence of values, we deduce that, in fact, the inequality in (\ref{eq:almosteverywhere}) holds for almost all $(u_1,u_2,\ldots,u_m)\in {\mathbb R}^m$. But then, using (\ref{eq:Fouriertransform}) and Theorem \ref{theorem:characterisation}, we see that  provided $m\ge 4s+1$, for almost all $(u_1,u_2,\ldots,u_m)\in {\mathbb R}^m$ there are $f_1,f_2,\ldots,f_m\in L^2({\mathbb R})$ such that 
\[f=\sum_{j=1}^m\left[
\left(
e^{ib_j   \left(\frac{\alpha-\beta}{2}\right)  } +e^{-ib_j\left(\frac{\alpha-\beta}{2}\right)}
\right)
\delta_0-\left(
e^{ib_j   \left(\frac{\alpha+\beta}{2}\right)  } \delta_{b_j}+e^{-ib_j\left(\frac{\alpha+\beta}{2}\right)}\delta_{-b_j}
\right)
\right]^s\ast f_j.\]
We deduce that (i) implies (ii) in Theorem \ref{theorem:main} and, by taking $m=4s+1$, we see that (i) implies (iii). 

We have now proved that  (i), (ii) and (iii) are equivalent. Also, the statement that (iii) is possible for almost all $(u_1,u_2,\ldots,u_{4s+1})\in {\mathbb R}^{4s+1}$ has been proved.

Finally, put $T=(D^2-i(\alpha+\beta)D-\alpha\beta I)^s$.  Then, if $g\in L^2({\mathbb R})$ we have $ (T(g))^ {\widehat {\ }}(x)=\hfill\break(-1)^s(x-\alpha)^s(x-\beta)^s{\widehat g}(x)$. Consequently, $\int_{-\infty}^{\infty}(x-\alpha)^{-2s}(x-\beta)^{-2s}| (T(g))^ {\widehat {\ }}(x)|^2\,dx<\infty$. As the multiplier of $T$ is $(-1)^s(x-\alpha)^{-2s}(x-\beta)^{-2s}$, it is easy to see that there is $K>0$ such that
\[||T(g)||_{\alpha,\beta,s}^2=\langle T(g),T(g)\rangle_{\alpha,\beta,s} \le K  \int_{-\infty}^{\infty} (1+x^{4s})|{\widehat g}(x)|^2\,dx\le K||g||_{{\mathbb R},2s}^2<\infty,\]
and it follows that $T$ is bounded from $W^{2s}({\mathbb R})$ into ${\mathcal D}_{\alpha,\beta,s}({\mathbb R})$. As the multiplier of $T$ vanishes only at the two points $\alpha$ and $\beta$, $T$ is injective on $W^{2s}({\mathbb R})$. Finally, if $h\in L^2 ({\mathbb R})$ is such that $\int_{-\infty}^{\infty}(x-\alpha)^{-2s}(x-\beta)^{-2s}| {\widehat h}(x)|^2\,dx<\infty$,  we may let   $g\in L^2({\mathbb R})$ be the function such that ${\widehat g}(x)=(-1)^s(x-\alpha)^{-s}(x-\beta)^{-s}{\widehat h}(x).$ It is easy to see that $g\in W^{2s}({\mathbb R})$ and that $T(g)=h$. Consequently, $T$ maps $W^{2s}({\mathbb R})$ onto  ${\mathcal D}_{\alpha,\beta,s}({\mathbb R})$, and it follows that $T$ is a bounded invertible linear operator from $W^{2s}({\mathbb R})$ onto ${\mathcal D}_{\alpha,\beta,s}({\mathbb R})$. This completes the proof of Theorem \ref{theorem:main}.\hfill$\square$

Note that  an alternative proof of Theorem \ref{theorem:main} for the special case  $\alpha=\beta$ may be derived from the identity (\ref{eq:characterisation2}), which was proved originally in \cite{nillsen1} and \cite{nillsen2}.    In \cite{meisters1} Meisters and Schmidt   showed that every translation invariant linear form on   $L^2({\mathbb T})$ is continuous, but in \cite{meisters2}  it was shown that there are discontinuous translation invariant linear forms on $L^2({\mathbb R})$, and this latter result may also be deduced from  the subsequent identity (\ref{eq:characterisation2}). 
 
 {\bf Definition.}   Let $\alpha,\beta\in {\mathbb R}$ and let  $s\in {\mathbb N}$. Then a linear form $T$ on $L^2(G)$ is called $(\alpha, \beta,s)$-\emph{invariant} if, for all $f\in L^2({\mathbb R})$  and $u\in L^2({\mathbb R})$, 
\[T\left(\,\left[
\left(e^{iu   \left(\frac{\alpha-\beta}{2}\right)  } +e^{-iu\left(\frac{\alpha-\beta}{2}\right)}
\right)\delta_0
-\left(
e^{iu   \left(\frac{\alpha+\beta}{2}\right)  } \delta_{u}+e^{-iu_j\left(\frac{\alpha+\beta}{2}\right)}\delta_{-u_j}
\right)
\right]^s\,\ast f\right)=T(f).\]
When $\alpha,\beta\in {\mathbb Z}$, we may also introduce the notion of $(\alpha,\beta,s)$-invariant linear forms on $L^2({\mathbb T})$.  It was shown in \cite{nillsen3} that an $(\alpha,\beta,1)$-invariant linear form on $L^2({\mathbb T})$ is continuous and, in fact,  any $(\alpha,\beta,s)$-invariant linear form on $L^2({\mathbb T})$ is continuous (proved using the same technique as in \cite{nillsen3} for the case $s=1$). Together with the preceding comments, the  following corollary to Theorem \ref{theorem:main} shows that the situation pertaining to translation invariant linear forms on $L^2({\mathbb T})$ and $L^2({\mathbb R})$ is mirrored by that for $(\alpha,\beta,s)$-invariant linear forms on $L^2({\mathbb T})$  and $L^2({\mathbb R})$.
\begin{corollary} Let $\alpha,\beta\in {\mathbb R}$ and let $s\in {\mathbb N}$. Then, there are discontinuous $(\alpha,\beta,s)$-invariant linear forms  on $L^2({\mathbb R})$.
\end{corollary}
{\bf Proof.} We see from the definitions that if $T$ is a linear form on $L^2({\mathbb R})$, then $T$ is $(\alpha,\beta,s)$-invariant if and only if $T$ vanishes on ${\mathcal D}_{\alpha,\beta,s}({\mathbb R})$. However, it is  consequence of Theorem \ref{theorem:main} that ${\mathcal D}_{\alpha,\beta,s}({\mathbb R})$ has infinite algebraic codimension in $L^2({\mathbb R})$. Consequently there are  discontinuous linear forms on $L^2({\mathbb R})$ that vanish on ${\mathcal D}_{\alpha,\beta,s}({\mathbb R})$, and such forms are also $(\alpha,\beta,s)$ invariant.\hfill$\square$
  
\end{document}